\documentclass{article}
\usepackage{amsfonts}
\usepackage{amsmath}

\newtheorem{theorem}{Theorem}

\newtheorem{remark}[theorem]{Remark}

\input{tcilatex}
\begin{document}

\title{Riemann-Lagrange geometry for starfish/coral dynamical system}
\author{Mircea Neagu}
\date{}
\maketitle

\begin{abstract}
In this paper we develop the Riemann-Lagrange geometry, in the sense of
nonlinear connection, d-torsions, d-curvatures and jet Yang-Mills entity,
associated with the dynamical system concerning social interaction in
colonial organisms.
\end{abstract}



\section{Social interactions in colonial organisms}

Let $m\geq 2$ be an integer. We introduce social interactions for
starfish/coral dynamics as follows \cite{Ant-Auger-Brad}:%
\begin{equation}
\left\{ _{{}}%
\begin{array}{lll}
\dfrac{dN^{1}}{dt} & = & \lambda _{1}N^{1}-\alpha _{1}\left( N^{1}\right)
^{2}-\alpha _{2}\left( \dfrac{m}{m-1}\right) \cdot N^{1}N^{2}+\medskip \\ 
& + & \dfrac{\alpha _{1}}{m-1}\left( \dfrac{N^{2}}{N^{1}}\right) ^{m-2}\cdot
\left( N^{2}\right) ^{2}-\delta _{1}FN^{1}\medskip \\ 
\dfrac{dN^{2}}{dt} & = & \lambda _{2}N^{2}-\alpha _{2}\left( N^{2}\right)
^{2}-\alpha _{1}\left( \dfrac{m}{m-1}\right) \cdot N^{1}N^{2}+\medskip \\ 
& + & \dfrac{\alpha _{2}}{m-1}\left( \dfrac{N^{1}}{N^{2}}\right) ^{m-2}\cdot
\left( N^{1}\right) ^{2}-\delta _{2}FN^{2}\medskip \\ 
\dfrac{dF}{dt} & = & \beta F\left( N^{1}+N^{2}\right) +\gamma F^{2}-\rho F,%
\end{array}%
\right.  \label{system}
\end{equation}%
where

\begin{itemize}
\item $\alpha _{1},$ $\alpha _{2},$ $\lambda _{1},$ $\lambda _{2},$ $\delta
_{1},$ $\delta _{2},$ $\beta ,$ $\gamma ,$ $\rho $ are positive coefficients;

\item $N^{1},$ $N^{2}$ are coral densities;

\item $F$ is the starfish density;

\item $\lambda _{1}$ and $\lambda _{2}$ are growth rates;

\item $\lambda _{1}/\alpha _{1}$ and $\lambda _{2}/\alpha _{2}$ are single
species carrying capacities;

\item $\beta ,$ $\delta _{1}$ and $\delta _{2}$ are the interaction
coefficients for starfish preying on corals;

\item $\gamma $ is the coefficient of starfish aggregation.
\end{itemize}

Note that $m$ is the effect of increasing the social parameter. If we set $%
m=2$, we obtain the (2 corals/1 starfish)-model of Antonelli and Kazarinoff 
\cite{Ant-Kaz}, in which every term of degree greater than one is quadratic.
It is $m\geq 3$ which forces the social interaction terms to be nonquadratic.

The dynamical system (\ref{system}) can be extended to a dynamical system of
order two coming from a first order Lagrangian of least square type. This
extension is called in the literature of specialty as the geometric
dynamical system (see Udri\c{s}te \cite{Udriste}).

\section{The Riemann-Lagrange geometry}

The system (\ref{system}) can be regarded on the 1-jet space $J^{1}(\mathbb{R%
},\mathbb{R}^{3})$, whose coordinates are%
\begin{equation*}
\left( t,\text{ }x^{1}=N^{1},\text{ }x^{2}=N^{2},\text{ }x_{3}=F,\text{ }%
y_{1}^{1}=\dfrac{dN^{1}}{dt},\text{ }y_{1}^{2}=\dfrac{dN^{2}}{dt},\text{ }%
y_{1}^{3}=\dfrac{dF}{dt}\right) .
\end{equation*}

\begin{remark}
We recall that the transformations of coordinates on the 1-jet space $J^{1}(%
\mathbb{R},\mathbb{R}^{3})$ are given by%
\begin{equation}
\left\{ 
\begin{array}{l}
\widetilde{t}=\widetilde{t}(t)\medskip \\ 
\widetilde{x}^{i}=\widetilde{x}^{i}(x^{j})\medskip \\ 
\widetilde{y}_{1}^{i}=\dfrac{\partial \widetilde{x}^{i}}{\partial x^{j}}%
\dfrac{dt}{d\widetilde{t}}\cdot y_{1}^{j},%
\end{array}%
\right.  \label{tr-of-coordinates}
\end{equation}%
where $i,j=\overline{1,3}.$
\end{remark}

In this context, the solutions of class $C^{2}$ of the system (\ref{system})
are the global minimum points of the jet least square Lagrangian \cite%
{Balan-Neagu}, \cite{Ne-Udr}%
\begin{equation}
L=\left( y_{1}^{1}-X_{(1)}^{(1)}\left( N^{1},N^{2},F\right) \right)
^{2}+\left( y_{1}^{2}-X_{(1)}^{(2)}\left( N^{1},N^{2},F\right) \right) ^{2}+
\label{Lagrangian}
\end{equation}%
\begin{equation*}
+\left( y_{1}^{3}-X_{(1)}^{(3)}\left( N^{1},N^{2},F\right) \right) ^{2},
\end{equation*}%
where%
\begin{equation*}
\begin{array}{lll}
X_{(1)}^{(1)}\left( N^{1},N^{2},F\right) & = & \lambda _{1}N^{1}-\alpha
_{1}\left( N^{1}\right) ^{2}-\alpha _{2}\left( \dfrac{m}{m-1}\right) \cdot
N^{1}N^{2}+\medskip \\ 
& + & \dfrac{\alpha _{1}}{m-1}\left( \dfrac{N^{2}}{N^{1}}\right) ^{m-2}\cdot
\left( N^{2}\right) ^{2}-\delta _{1}FN^{1}\medskip%
\end{array}%
\end{equation*}%
\begin{equation*}
\begin{array}{lll}
X_{(1)}^{(2)}\left( N^{1},N^{2},F\right) & = & \lambda _{2}N^{2}-\alpha
_{2}\left( N^{2}\right) ^{2}-\alpha _{1}\left( \dfrac{m}{m-1}\right) \cdot
N^{1}N^{2}+\medskip \\ 
& + & \dfrac{\alpha _{2}}{m-1}\left( \dfrac{N^{1}}{N^{2}}\right) ^{m-2}\cdot
\left( N^{1}\right) ^{2}-\delta _{2}FN^{2}\medskip \\ 
X_{(1)}^{(3)}\left( N^{1},N^{2},F\right) & = & \beta F\left(
N^{1}+N^{2}\right) +\gamma F^{2}-\rho F,%
\end{array}%
\end{equation*}

\begin{remark}
The solutions of class $C^{2}$ of the system (\ref{system}) are solutions of
the Euler-Lagrange equations attached to the the jet least square Lagrangian
(\ref{Lagrangian}), namely (geometric dynamics)%
\begin{equation*}
\frac{\partial L}{\partial x^{i}}-\frac{d}{dt}\left( \frac{\partial L}{%
\partial y_{1}^{i}}\right) =0,\text{ }\forall \text{ }i=\overline{1,3}.
\end{equation*}
\end{remark}

But, the jet least square Lagrangian (\ref{Lagrangian}) provides us with an
entire Riemann-Lagrange geometry on the 1-jet space $J^{1}(\mathbb{R},%
\mathbb{R}^{2})$, in the sense of nonlinear connection, d-torsions,
d-curvatures and jet Yang-Mills entity. Let us use the notation%
\begin{equation*}
J(X_{(1)})=\left( \frac{\partial X_{(1)}^{(i)}}{\partial x^{j}}\right) _{i,j=%
\overline{1,3}}=\left( 
\begin{array}{ccc}
J_{11} & J_{12} & J_{13}\medskip \\ 
J_{21} & J_{22} & J_{23}\medskip \\ 
J_{31} & J_{32} & J_{33}%
\end{array}%
\right) ,
\end{equation*}%
where\medskip

$J_{11}=\lambda _{1}-2\alpha _{1}N^{1}-\alpha _{2}\left( \dfrac{m}{m-1}%
\right) \cdot N^{2}-\alpha _{1}\left( \dfrac{m-2}{m-1}\right) \dfrac{\left(
N^{2}\right) ^{m}}{\left( N^{1}\right) ^{m-1}}-\delta _{1}F,\medskip$

$J_{12}=-\alpha _{2}\left( \dfrac{m}{m-1}\right) \cdot N^{1}+\alpha
_{1}\left( \dfrac{m}{m-1}\right) \dfrac{\left( N^{2}\right) ^{m-1}}{\left(
N^{1}\right) ^{m-2}},\medskip $

$J_{13}=-\delta _{1}N^{1},\medskip $

$J_{21}=-\alpha _{1}\left( \dfrac{m}{m-1}\right) \cdot N^{2}+\alpha
_{2}\left( \dfrac{m}{m-1}\right) \dfrac{\left( N^{1}\right) ^{m-1}}{\left(
N^{2}\right) ^{m-2}},\medskip $

$J_{22}=\lambda _{2}-2\alpha _{2}N^{2}-\alpha _{1}\left( \dfrac{m}{m-1}%
\right) \cdot N^{1}-\alpha _{2}\left( \dfrac{m-2}{m-1}\right) \dfrac{\left(
N^{1}\right) ^{m}}{\left( N^{2}\right) ^{m-1}}-\delta _{2}F,\medskip$

$J_{23}=-\delta _{2}N^{2},\medskip $

$J_{31}=\beta F,\medskip $

$J_{32}=\beta F,\medskip $

$J_{33}=\beta \left( N^{1}+N^{2}\right) +2\gamma F-\rho .\medskip $

Following the geometrical ideas from Balan and Neagu \cite{Balan-Neagu}, we
obtain the following geometrical results:

\begin{theorem}
\emph{(i)} The canonical nonlinear connection on $J^{1}(\mathbb{R},\mathbb{R}%
^{3})$, produced by the system (\ref{system}), has the local components%
\begin{equation*}
\Gamma =\left( M_{(1)1}^{(i)}=0,\text{ }N_{(1)j}^{(i)}\right) ,
\end{equation*}%
where $N_{(1)j}^{(i)}$ are the entries of the skew-symmetric matrix%
\begin{equation*}
N_{(1)}=-\frac{1}{2}\left[ J(X_{(1)})-\text{ }^{T}J(X_{(1)})\right] =\left( 
\begin{array}{ccc}
N_{(1)1}^{(1)} & N_{(1)2}^{(1)} & N_{(1)3}^{(1)}\medskip \\ 
N_{(1)1}^{(2)} & N_{(1)2}^{(2)} & N_{(1)3}^{(2)}\medskip \\ 
N_{(1)1}^{(3)} & N_{(1)2}^{(3)} & N_{(1)3}^{(3)}%
\end{array}%
\right) ,
\end{equation*}%
where\medskip

$N_{(1)1}^{(1)}=N_{(1)2}^{(2)}=N_{(1)3}^{(3)}=0,\medskip $

$N_{(1)2}^{(1)}=-N_{(1)1}^{(2)}=\left( \dfrac{m}{m-1}\right) \left( \alpha
_{1}N^{2}-\alpha _{2}N^{1}\right) +\medskip $

\hspace{25mm}$+\left( \dfrac{m}{m-1}\right) \left[ \alpha _{2}\dfrac{\left(
N^{1}\right) ^{m-1}}{\left( N^{2}\right) ^{m-2}}-\alpha _{1}\dfrac{\left(
N^{2}\right) ^{m-1}}{\left( N^{1}\right) ^{m-2}}\right] ,\medskip $

$N_{(1)3}^{(1)}=-N_{(1)1}^{(3)}=\beta F+\delta _{1}N^{1},\medskip _{{}}$

$N_{(1)3}^{(2)}=-N_{(1)2}^{(3)}=\beta F+\delta _{2}N^{2}.\medskip $

\emph{(ii)} All adapted components of the canonical generalized Cartan
connection $C\Gamma $, produced by the system (\ref{system}), are zero.

\emph{(iii)} The effective adapted components $R_{(1)jk}^{(i)}$ of the
torsion d-tensor $\mathbf{T}$ of the canonical generalized Cartan connection 
$C\Gamma $, produced by the system (\ref{system}), are the entries of the
following skew-symmetric matrices:%
\begin{equation*}
R_{(1)1}=\left( R_{(1)j1}^{(i)}\right) _{i,j=\overline{1,3}}=\frac{\partial
N_{(1)}}{\partial N^{1}}=\left( 
\begin{array}{ccc}
0 & \dfrac{\partial N_{(1)2}^{(1)}}{\partial N^{1}} & \delta _{1}\medskip \\ 
-\dfrac{\partial N_{(1)2}^{(1)}}{\partial N^{1}} & 0 & 0\medskip \\ 
-\delta _{1} & 0 & 0%
\end{array}%
\right) ,
\end{equation*}%
where\medskip

\noindent$\dfrac{\partial N_{(1)2}^{(1)}}{\partial N^{1}}=\left( \dfrac{m}{%
m-1}\right) \left[ -\alpha _{2}+\alpha _{2}\left( m-1\right) \left( \dfrac{%
N^{1}}{N^{2}}\right) ^{m-2}+\alpha _{1}\left( m-2\right) \left( \dfrac{N^{2}%
}{N^{1}}\right) ^{m-1}\right] ;\medskip$%
\begin{equation*}
R_{(1)2}=\left( R_{(1)j2}^{(i)}\right) _{i,j=\overline{1,3}}=\frac{\partial
N_{(1)}}{\partial N^{2}}=\left( 
\begin{array}{ccc}
0 & \dfrac{\partial N_{(1)2}^{(1)}}{\partial N^{2}} & 0\medskip \\ 
-\dfrac{\partial N_{(1)2}^{(1)}}{\partial N^{2}} & 0 & \delta _{2}\medskip
\\ 
0 & -\delta _{2} & 0%
\end{array}%
\right) ,
\end{equation*}%
where\medskip

\noindent$\dfrac{\partial N_{(1)2}^{(1)}}{\partial N^{2}}=\left( \dfrac{m}{%
m-1}\right) \left[ \alpha _{1}-\alpha _{2}\left( m-2\right) \left( \dfrac{%
N^{1}}{N^{2}}\right) ^{m-1}-\alpha _{1}\left( m-1\right) \left( \dfrac{N^{2}%
}{N^{1}}\right) ^{m-2}\right] ;\medskip$

\begin{equation*}
R_{(1)3}=\left( R_{(1)j3}^{(i)}\right) _{i,j=\overline{1,3}}=\frac{\partial
N_{(1)}}{\partial F}=\left( 
\begin{array}{ccc}
0 & 0 & \beta \medskip \\ 
0 & 0 & \beta \medskip \\ 
-\beta & -\beta & 0%
\end{array}%
\right) .
\end{equation*}

\emph{(iv)} All adapted components of the curvature d-tensor $\mathbf{R}$ of
the canonical generalized Cartan connection $C\Gamma $, produced by the
system (\ref{system}), cancel.

\emph{(v)} The geometric electromagnetic-like distinguished $2$-form,
produced by the system (\ref{system}), is given by%
\begin{equation*}
\mathbb{F}=F_{(i)j}^{(1)}\delta y_{1}^{i}\wedge dx^{j},
\end{equation*}%
where%
\begin{equation*}
\delta y_{1}^{i}=dy_{1}^{i}-N_{(1)j}^{(i)}dx^{j},\quad \forall \text{ }i=%
\overline{1,3},
\end{equation*}%
and the adapted components $F_{(i)j}^{(1)}$ are the entries of the
skew-symmetric matrix%
\begin{equation*}
F^{(1)}=\left( F_{(i)j}^{(1)}\right) _{i,j=\overline{1,3}}=-N_{(1)}.
\end{equation*}

\emph{(vi)} The jet geometric Yang-Mills entity, produced by the system (\ref%
{system}), is given by the formula%
\begin{equation*}
\mathcal{EYM}(t)=\left[ F_{(1)2}^{(1)}\right] ^{2}+\left[ F_{(1)3}^{(1)}%
\right] ^{2}+\left[ F_{(2)3}^{(1)}\right] ^{2}=
\end{equation*}%
\begin{equation*}
=\left\{ \left( \dfrac{m}{m-1}\right) \left( \alpha _{1}N^{2}-\alpha
_{2}N^{1}\right) +\left( \dfrac{m}{m-1}\right) \left[ \alpha _{2}\dfrac{%
\left( N^{1}\right) ^{m-1}}{\left( N^{2}\right) ^{m-2}}-\alpha _{1}\dfrac{%
\left( N^{2}\right) ^{m-1}}{\left( N^{1}\right) ^{m-2}}\right] \right\} ^{2}+
\end{equation*}%
\begin{equation*}
+\left( \beta F+\delta _{1}N^{1}\right) ^{2}+\left( \beta F+\delta
_{2}N^{2}\right) ^{2}.
\end{equation*}
\end{theorem}

\begin{remark}
We recall that, under a transformation of coordinates (\ref%
{tr-of-coordinates}), the local components of the nonlinear connection obey
the rules (see \cite{Balan-Neagu}, ,\cite{Buc-Mir}, \cite{Mir-An})%
\begin{equation*}
\widetilde{M}_{(1)1}^{(k)}=M_{(1)1}^{(j)}\left( \frac{dt}{d\widetilde{t}}%
\right) ^{2}\dfrac{\partial \widetilde{x}^{k}}{\partial x^{j}}-\frac{dt}{d%
\widetilde{t}}\frac{\partial \widetilde{y}_{1}^{k}}{\partial t},
\end{equation*}%
\begin{equation*}
\widetilde{N}_{(1)l}^{(k)}=N_{(1)i}^{(j)}\frac{dt}{d\widetilde{t}}\dfrac{%
\partial x^{i}}{\partial \widetilde{x}^{l}}\dfrac{\partial \widetilde{x}^{k}%
}{\partial x^{j}}-\dfrac{\partial x^{i}}{\partial \widetilde{x}^{l}}\frac{%
\partial \widetilde{y}_{1}^{k}}{\partial x^{i}}.
\end{equation*}
\end{remark}

\noindent \textbf{Open problem.} There exist real meanings in trophodynamics
for the geometrical objects constructed in this paper?

Mircea Neagu

Department of Mathematics and Informatics

Transilvania University of Bra\c{s}ov

Blvd. Iuliu Maniu, No. 50, Bra\c{s}ov 500091, Romania

email: \textit{mircea.neagu@unitbv.ro}

\end{document}